# Modeling dynamic crack branching in unsaturated porous media through multi-phase micro-periporomechanics

Hossein Pashazad[a], Xiaoyu Song[a,*]

[a]*Engineering School of Sustainable Infrastructure and Environment*
*University of Florida, Gainesville, FL, USA*

**Abstract**

Dynamic crack branching in unsaturated porous media holds significant relevance in various fields, including geotechnical engineering, geosciences, and petroleum engineering. This article presents a numerical investigation into dynamic crack branching in unsaturated porous media using a recently developed coupled micro-periporomechanics paradigm. This paradigm extends the periporomechanics model by incorporating the micro-rotation of the solid skeleton. Within this framework, each material point is equipped with three degrees of freedom: displacement, micro-rotation, and fluid pressure. Consistent with the Cosserat continuum theory, a length scale associated with the micro-rotation of material points is inherently integrated into the model. This study encompasses several key aspects: (1) Validation of the coupled micro-periporomechanics paradigm for effectively modeling crack branching in deformable porous media, (2) Examination of the transition from a single branch to multiple branches in porous media under drained conditions, (3) Simulation of single crack branching in unsaturated porous media under dynamic loading conditions, and (4) Investigation of multiple crack branching in unsaturated porous media under dynamic loading conditions. The numerical results obtained in this study are systematically analyzed to elucidate the factors that influence dynamic crack branching in porous media subjected to dynamic loading. Furthermore, the comprehensive numerical findings underscore the efficacy and robustness of the coupled micro-periporomechanics paradigm in accurately modeling dynamic crack branching in variably saturated porous media.

## 1. Introduction

Dynamic crack branching in unsaturated porous media holds significant relevance in various fields, including geotechnical engineering, geosciences, and petroleum engineering [1–9]. Cracking in unsaturated soils can compromise the structural integrity of infrastructure built upon unsaturated soils. This cracking can be induced by matric suction in unsaturated soils, leading to volume shrinkage and the formation of tensile cracks [10, 11]. Such cracks, in turn, directly impact the soil's bearing capacity, thereby affecting the stability of foundations for structures situated on such soil. In mountainous regions, landslides are common, encompassing various modes such as creeping motion, initial failures, rapid sliding, and transitions to very rapid movement. Surface cracks can trigger these landslides [12–16]. In the field of gas and soil production engineering, hydraulic fracturing has emerged as a crucial stimulation technique to enhance the creation of a conductive network of fractures, thereby boosting the production of unconventional natural resources [17–20]. This process involves injecting high-pressure fluids into the bedrock formation to either widen existing fractures or generate new ones. This operation necessitates a close interaction between the solid framework and the flow of fluids within interconnected voids, resulting in a highly complex multiphase interplay. Under specific conditions, the fluid-induced fracture may change its course or bifurcate into multiple branches [21, 22]. Although experimental and numerical studies have provided evidence of crack branching during hydraulic fracturing, a consensus on the precise factors influencing this phenomenon is yet to be reached [22, 23]. Consequently, the accurate prediction of hydraulic fracturing, particularly when considering inertial effects, remains an ongoing research area that warrants further exploration. In this article, as a new contribution, we numerically investigate the dynamic crack branching in unsaturated porous media through a recently formulated coupled micro-periporomechanics ($\mu$PPM) paradigm [24, 25] that extends the periporomechanics (PPM) model by considering the micro-rotation of the solid skeleton. Next, we briefly review the physical experiment study of crack branching in porous media.

---

*Corresponding author
Email address:* xysong@ufl.edu (Xiaoyu Song)

Crack branching in porous media, such as clay, has correlated with the clay layers' mechanical properties [26]. When the stress within a fracture zone exceeds a critical threshold in porous materials, the material cannot dissipate the energy efficiently, leading to crack branching, which occurs when a small critical energy release rate is surpassed. The phenomenon of crack branching is particularly notable when the crack reaches a critical speed of propagation [27]. It is important to note that the failure mode can transition from mode I to mixed modes with an increase in loading rate. At high crack propagation speeds, inertia forces at the crack tip impede crack propagation, resulting in branching (e.g., [27]). Experimental studies have consistently demonstrated that crack growth velocity in porous materials is lower than the Rayleigh wave velocity (e.g., [28, 29]). The crack branching criteria necessitate a critical dynamic stress intensity factor and consideration of the crack's curvature [30]. These criteria are valuable tools for predicting crack branching in dynamic brittle fracture tests [30]. In an experimental study [31], crack bifurcation, where a crack extends into multiple branches, was observed once the critical velocity leading to the initiation of velocity oscillations is exceeded. These experiments were conducted on delicate, nearly two-dimensional layers of brittle material. The instability analysis outcomes provided insight into the fracture process, emphasizing the presence and progressive development of instability as a precursor to intricate microscopic branching [31]. Further research on the crack branching in shale under tensile stress is presented in [32], where a curved specimen containing a deliberately created artificial notch along its curved edge was employed to observe the sequence of damage and crack propagation during brittle fracturing of shale. Notably, crack growth often entailed the cessation or closure of former branch cracks due to elastic recovery and induced compressive stress [32]. While physical testing is essential for studying crack branching in porous media, novel numerical modeling is equally vital in probing dynamic crack branching. Next, we present a brief review of the $\mu$PPM, which will be used for modeling dynamic crack branching in this study. For other numerical methods for modeling crack branching in porous media, we refer to the literature (e.g., [33–35], among others).

PPM is a nonlocal formulation of classical poromechanics [36–38] in the form of integro-differential equations through peridynamic states and the effective force concept [39, 40]. In PPM, the porous media is postulated to consist of a finite number of mixed material points that, within a finite distance called horizon, have direct poromechanical interaction. The PPM paradigm has been numerically implemented through the total and updated Lagrangian meshfree method in space and the monolithic/fractional-step implicit and explicit Newmark schemes in time. For a comparison between PPM and other numerical methods for modeling porous media, we refer to [41]. Within the PPM framework, multiphase discontinuities can naturally emerge based on field equations and material models [40]. Classical constitutive models for porous media can be incorporated into the PPM framework using the stabilized multiphase constitutive correspondence principle [40, 42]. The computational meshfree PPM method has been used to study instability, large deformation, and fracturing in variably saturated porous media under static and dynamic loads [42–50]. The $\mu$PPM paradigm has been recently developed to extend the original PPM paradigm by incorporating the micro-rotation of the solid skeleton following the Cosserat continuum theory [51–53]. In $\mu$PPM, each material point is equipped with three degrees of freedom, i.e., displacement, micro-rotation, and fluid pressures. Consistent with the Cosserat continuum theory, a micro-structure-based length scale associated with the micro-rotation of material points is inherently integrated into the $\mu$PPM paradigm. The stabilized micro-polar multiphase micro-polar constitutive correspondence principle has been formulated to incorporate the classical micro-polar material models for porous media into the new $\mu$PPM paradigm. In [25], we have numerically implemented the $\mu$PPM paradigm through the hybrid Lagrangian-Eulerian meshfree method and an explicit-explicit fractional-step algorithm.

In this study, we, for the first time, utilize the newly formulated coupled $\mu$PPM paradigm [24, 25] to numerically investigate the dynamic crack branching in porous media accounting for the micro-structure of porous materials (i.e., micro-rotations). In this journey, we first validate the coupled $\mu$PPM for modeling crack branching in deformable porous media. Second, we study the crack branching from a single branch to multiple branches in dry porous media. Third, we simulate the single crack branching in unsaturated porous media under dynamic loading. Lastly, we investigate the multiple crack branching in unsaturated porous media under dynamic loading. The numerical results are analyzed to show the factors influencing dynamic crack branching in porous media under dynamic loading. Furthermore, our comprehensive numerical results have demonstrated the efficacy and robustness of the coupled $\mu$PPM paradigm in modeling dynamic crack branching in unsaturated porous media. As implied by our numerical results, we note that the Cosserat length scale contributes to alleviating the dispersion issue with the standard PD models for solids.

The remainder of this article is organized as follows. Section 2 presents the mathematical formulation of the coupled $\mu$PPM paradigm and its numerical implementation. Section 3 presents numerical examples to validate the coupled mesh-free micro-PPM paradigm and demonstrate its efficacy and robustness in modeling crack branching in unsaturated porous media under dynamic loads, followed by a summary of the present study in Section 4. For the sign convention, the assumption in continuum mechanics [54] is adopted, i.e., the tensile force and deformation under tension are positive. For pore fluid, compression is positive, and tension is negative.

## 2. Mathematical formulation

This section presents the mathematical formulation and numerical implementation of the mesh-free fracturing $\mu$PPM paradigm for modeling the dynamic crack branching in unsaturated porous media.

### 2.1. Governing equations for the fracturing μPPM paradigm

In $\mu$PPM, the unsaturated porous material is represented by a collection of mixed material points. The material points at a finite distance called horizon $\delta$ have direct poromechanical interactions. This study assumes that the material point interacts with all material points in a spherical domain $H$ centered at the material point with a radius of $\delta$. The mixed material points have three types of degrees of freedom, i.e., displacement, micro-rotation, and pore fluid pressures. It is assumed that the solid skeleton has micro-rotations and the fluid phase has no micro-rotation, i.e., non-polar. In this study, the unsaturated porous media is assumed to comprise three phases, i.e., solid skeleton, pore water, and pore air. By assuming a weightless pore air phase, the density of unsaturated porous media $\rho$ is written as

$$\rho = (1-\phi)\rho_s + S_r \phi \rho_w, \tag{1}$$

where $\rho_s$ is the intrinsic mass density of solid phase, $\rho_w$ in the intrinsic mass density of water, $\varphi$ is the porosity, and $S_r$ is the degree of saturation. By assuming passive air pressure (i.e., atmospheric air pressure), the matric suction $s$ is defined as negative pore pressure under unsaturated conditions. Following the classical coupled unsaturated poromechanics in $\mu$PPM, we note that the solid skeleton is described using the Lagrangian coordinate system, and the fluid phase is described using the Eulerian coordinate system relative to the solid skeleton.

Next, we introduce the kinematics of two material points in the $\mu$PPM paradigm. Figure 1 plots the kinematics of two mixed material points in an unsaturated porous material body in the $\mu$PPM paradigm. Let $\boldsymbol{x}$ and $\boldsymbol{x}'$ represent two mixed material points in the reference configuration of a porous body. The micro-polar poromechanics bond between the two points is defined as $\boldsymbol{\xi} = \boldsymbol{x}' - \boldsymbol{x}$. The deformation vector state and the displacement vector state on the bond $\boldsymbol{\xi}$ are defined as

$$\underline{\mathscr{Y}} = \boldsymbol{y}' - \boldsymbol{y}, \tag{2}$$

$$\underline{\mathscr{U}} = \boldsymbol{u}' - \boldsymbol{u}, \tag{3}$$

where $\boldsymbol{y}$ and $\boldsymbol{y}'$ are the spatial locations of the two mixed material points in the current configuration, respectively, and $\boldsymbol{u}$ and $\boldsymbol{u}'$ are the displacements of the two material points, respectively.

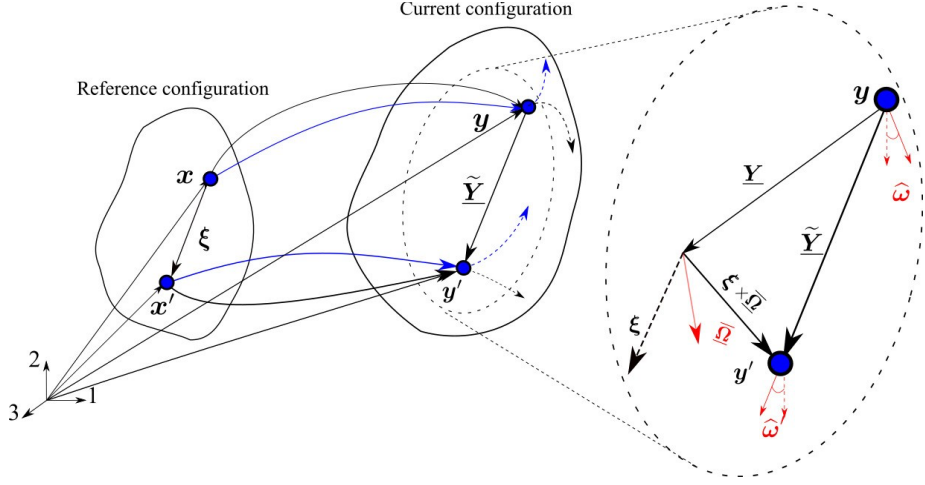

Figure 1: Kinematics of the mixed material points in the μPPM paradigm.

The micro-rotation vector state and the composite displacement vector state are defined as

$$\underline{\Omega} = \widehat{\omega}' - \widehat{\omega}, \quad (4)$$

$$\underline{\widetilde{\mathscr{U}}} = \underline{\mathscr{U}} - \overline{\Omega} \times \underline{\xi}. \quad (5)$$

where $\widehat{\omega}$ and $\widehat{\omega}'$ are the micro-rotations at $x$ and $x'$, respectively, and $\overline{\Omega}$ is the averaged micro-rotation vector state that reads

$$\overline{\Omega} = \frac{1}{2}\left(\widehat{\omega}' + \widehat{\omega}\right). \quad (6)$$

The fluid pressure state is defined as

$$\underline{\Phi} = p'_w - p_w, \quad (7)$$

where $p_w$ and $p'_w$ are the fluid pressures at the two material points $x$ and $x'$ in the current configuration, respectively. For notation simplicity, in the remaining presentation, a variable without a prime is associated with material point $x$, and a variable with a prime is associated with material point $x'$. Next, we present the governing equations for the fully coupled μPPM paradigm.

The governing equations consist of the motion equation, the moment balance equation, and the mass balance equation. The motion equation is written as

$$\rho \ddot{u} = \int_{\mathscr{H}} (\overline{\underline{\mathscr{T}}} - \overline{\underline{\mathscr{T}}}')d\mathscr{V}' - \int_{\mathscr{H}} (S_l \underline{\mathscr{T}}_l - S'_l \underline{\mathscr{T}}'_l)d\mathscr{V}' + \rho g, \quad (8)$$

where $\ddot{u}$ is acceleration, $g$ is gravitational acceleration, $\overline{\underline{T}}$ and $\overline{\underline{T}}'$ are the effective force states, $S_l$ and $S'_l$ are the degrees of saturation, and $\underline{T}_l$ and $\underline{T}'_l$ are the fluid force states. The term is expressed as

$$S_l \underline{\mathscr{T}}_l = \begin{cases} S_{r,f} \underline{\mathscr{T}}_f & \text{if } D > D_{cr} \ \& \ D' > D_{cr}, \\ S_r \underline{\mathscr{T}}_w, & \text{otherwise,} \end{cases} \quad (9)$$

where $\underline{T}_f$ is the fluid pressure state in a fractured point, $\underline{T}_w$ is the fluid pressure state in a bulk point, $S_{r,f}$ is the degree of saturation of a fractured point, $S_r$ is the degree of saturation of a bulk point, $D$ is the damage variable which is defined in the following section, and $D_{cr}$ is the critical damage variable. Note that the term $S'_l \underline{T}'_l$ at material point $x'$ in (8) can be expressed following (9). The degree of saturation can be determined through the soil-water retention curve [38]. Assuming passive air pressure, the soil-water retention curve is written as

$$S_r = \left[1 + \left(\frac{-p_w}{s_a}\right)^n\right]^{-m}, \tag{10}$$

where $s_a$, $m$, and $n$ are material parameters. The moment balance equation is written as

$$\mathcal{I}^s \ddot{\underline{\omega}} = \int_{\mathcal{H}} (\underline{\underline{M}} - \underline{\underline{M}}') \, d\mathcal{V}' + \frac{1}{2} \int_{\mathcal{H}} \underline{\mathcal{Y}} \times \left[ (\underline{\underline{T}} - S_l \underline{\underline{T}}_l) - (\underline{\underline{T}}' - S_l' \underline{\underline{T}}_l') \right] d\mathcal{V}' + \boldsymbol{l}, \tag{11}$$

where $\boldsymbol{I}^s$ the micro-inertia of the solid phase, $\ddot{\underline{\omega}}$ is the angular acceleration, $\underline{\underline{M}}$ and $\underline{\underline{M}}'$ are the moment states, and $\boldsymbol{l}$ is the body couple density.

In this study, it is assumed that the solid grain and water are incompressible. The mass balance equation in the bulk space is written as

$$\phi \frac{dS_r}{dt} + S_r \dot{\mathcal{V}} + \frac{1}{\rho_w} \int_{\mathcal{H}} (\underline{\mathcal{Q}} - \underline{\mathcal{Q}}') \, d\mathcal{V}' + Q_s = 0, \tag{12}$$

where $\dot{\mathcal{V}}$ is the volume change rate of the solid skeleton, $\underline{Q}$ and $\underline{Q}'$ are the fluid flow states, and $Q_s$ is a source term. It is noted that the micro-rotation of the solid phase does not affect the volume change rate of the solid phase [55]. Following (12), the mass balance equation in the fractured space can be written as

$$\frac{\partial S_{r,f}}{\partial t} + \frac{1}{\rho_w} \int_{\mathcal{H}} (\underline{\mathcal{Q}}_f - \underline{\mathcal{Q}}'_f) \, d\mathcal{V}' - Q_s = 0, \tag{13}$$

where $\underline{Q}_f$ and $\underline{Q}'_f$ are the fluid flow states of the material points in the fractured space. In (13), it is assumed that the porosity $\varphi = 1$ in the fractured space and the volume coupling term vanishes.

To determine $Q_s$, we assume that the direction of fluid flow from the bulk to the fractured space is normal to the fracture surface. Then, following the generalized Darcy's law, $Q_s$ can be written as

$$Q_s = 2 \left[ -\frac{k^r k_w}{\mu_w} \left( \frac{p_f - p_w}{d} \right) \right] / d, \tag{14}$$

where $p_w$ is the water pressure in bulk, $p_f$ is the water pressure in the fractured space, and $d$ is the edge dimension of a cubic material point in a uniform grid [47]. We note that for a material point and its neighbor points in the bulk, (8) and (16) degenerate into the following equations.

$$\rho \ddot{\boldsymbol{u}} = \int_{\mathcal{H}} (\underline{\underline{T}} - \underline{\underline{T}}') d\mathcal{V}' - \int_{\mathcal{H}} (S_r \underline{\underline{T}}_w - S_r' \underline{\underline{T}}'_w) d\mathcal{V}' + \rho \boldsymbol{g}, \tag{15}$$

$$\mathcal{I}^s \ddot{\underline{\omega}} = \int_{\mathcal{H}} (\underline{\underline{M}} - \underline{\underline{M}}') \, d\mathcal{V}' + \frac{1}{2} \int_{\mathcal{H}} \underline{\mathcal{Y}} \times \left[ (\underline{\underline{T}} - S_r \underline{\underline{T}}_w) - (\underline{\underline{T}}' - S_r' \underline{\underline{T}}'_w) \right] d\mathcal{V}' + \boldsymbol{l}. \tag{16}$$

In summary, the governing equations of the coupled μPPM consists of (8), (16), (12), and (13). To complete the mathematical framework, we adopt the classical micro-polar material models to determine the effective force, moment, and fluid flow states through the stabilized micro-polar multiphase constitutive correspondence principle in the following section.

*2.2. Stabilized multiphase μPPM correspondence principle*

This part presents the nonlocal constitutive models harnessing the classical micro-polar material models for the skeleton and the generalized non-polar Dacy's law for unsaturated fluid flow. In so doing, the stabilized multiphase μPPM correspondence principle [25] is used to determine the nonlocal strain tensor, wryness tensor, and fluid pressure gradient vector. Assuming the small deformation of the skeleton, the nonlocal strain tensor $\boldsymbol{\varepsilon}$ and the nonlocal wryness tensor $\boldsymbol{\kappa}$ can be written as

$$\varepsilon = \left[\int_{\mathcal{H}} \underline{\omega} \left(\widetilde{\underline{\mathcal{U}}} \otimes \underline{\xi}\right) d\mathcal{V}'\right] \mathcal{K}^{-1}, \tag{17}$$

$$\kappa = \left[\int_{H} \underline{\omega} \left(\underline{\Omega} \otimes \underline{\xi}\right) d\mathcal{V}'\right] \mathcal{K}^{-1}, \tag{18}$$

where $\underline{\omega}$ is the unit weighting function and $K$ is the shape tensor, defined as

$$\mathcal{K} = \int_{\mathcal{H}} \underline{\xi} \otimes \underline{\xi} d\mathcal{V}'. \tag{19}$$

Similarly, the nonlocal fluid pressure gradients in the bulk and the fractured space can be written as

$$\widetilde{\boldsymbol{\nabla}\Phi} = \left(\int_{\mathcal{H}} \omega \underline{\Phi}\underline{\xi} d\mathcal{V}'\right) \mathcal{K}^{-1}, \tag{20}$$

$$\widetilde{\boldsymbol{\nabla}\Phi_f} = \left(\int_{\mathcal{H}} \omega \underline{\Phi}_f \underline{\xi} d\mathcal{V}'\right) \mathcal{K}^{-1}, \tag{21}$$

where $\underline{\Phi}_f = p'_f - p_f$, and $p'_f$ and $p_f$ are the fluid pressures in the fractured material points.

Given (17) and (18), the effective stress tensor and the couple stress tensor can be computed from the classical micro-polar constitutive model for the solid skeleton. In this study, an elastic micro-polar elastic model is adopted as follows.

$$\overline{\sigma}_{ij} = \lambda \varepsilon_{kk} + (\mu + \mu_c)\varepsilon_{ij} + (\mu - \mu_c)\varepsilon_{ji}, \tag{22}$$

$$m_{ij} = \frac{1}{2}\mu l^2 \kappa_{ij}, \tag{23}$$

where $i, j, k = 1, 2, 3$, $\lambda$ is Lame's first elastic constant, $\mu$ is the shear modulus, $\mu_c$ is the micropolar shear modulus, and $l$ is the micropolar length scale [53]. Similarly, given (20) and (20), through the generalized Darcy's law, the unsaturated flow fluid flux vectors in the bulk and the fractured space can be written as

$$\boldsymbol{q} = -\frac{k^r k_w}{\mu_w}\widetilde{\boldsymbol{\nabla}\Phi}, \tag{24}$$

$$\boldsymbol{q}_f = -\frac{k^r_f k_f}{\mu_w}\widetilde{\boldsymbol{\nabla}\Phi_f}, \tag{25}$$

where $k_w$ and $k_r$ are the intrinsic and relative permeabilities of the bulk, respectively, $k_f$ and $k_{rf}$ are the intrinsic and relative permeabilities of the fractured space, respectively, and $\mu_w$ is the viscosity of water. The relative permeabilities [38] can be written as

$$k^r = S_r^{1/2}\left[1 - \left(1 - S_r^{1/m}\right)^m\right]^2, \tag{26}$$

$$k^r_f = S_{r,f}^{1/2}\left[1 - \left(1 - S_{r,f}^{1/m}\right)^m\right]^2, \tag{27}$$

where $m$ is the material parameter defined in (10). In this study, it is assumed the intrinsic permeability in the fractured space is determined through the cubic law [6] as

$$k_f = \frac{a_f^2}{12}, \tag{28}$$

where $a_f$ is the crack width at fracture space.

Given (22) and (23), the stabilized effective force and moment states [42] can be written as

$$\underline{\overline{\mathcal{T}}} = \underline{\omega}\overline{\sigma}\mathcal{K}^{-1}\underline{\xi} + \frac{\mathscr{GC}_1}{\omega_0}\underline{\omega}\underline{\mathscr{R}}_1, \tag{29}$$

$$\underline{\mathscr{M}} = \underline{\omega}m\mathcal{K}^{-1}\underline{\xi} + \frac{\mathscr{GC}_2}{\omega_0}\underline{\omega}\underline{\mathscr{R}}_2 \tag{30}$$

where $G$ is a positive number on the order of 1, $C_1$ and $C_2$ are two material parameters,

$$\omega_0 = \int_{\mathcal{H}} \omega \underline{\xi} d\mathcal{V}', \tag{31}$$

and $\underline{R}_1$ and $\underline{R}_2$ are the non-uniform composite displacement state and the non-uniform micro-rotation state, respectively. The latter two terms are defined as

$$\underline{\mathcal{R}}_1 = \widetilde{\underline{\mathcal{U}}} - \varepsilon \underline{\xi}, \tag{32}$$

$$\underline{\mathcal{R}}_2 = \underline{\Omega} - \kappa \underline{\xi}. \tag{33}$$

Under the three dimensions, $C_1$ and $C_2$ can be written as

$$\mathcal{C}_1 = \frac{12\mathcal{D}}{|\underline{\xi}|^3}, \tag{34}$$

$$\mathcal{C}_2 = \frac{\mathcal{D}}{|\underline{\xi}|}, \tag{35}$$

where $D$ is variable that depends on the horizon and material properties [42]. For a three-dimensional case, it can be determined by

$$\mathcal{D} = \frac{E(1 - 4\nu)}{4\pi\delta^2(1 - \nu - 2\nu^2)}, \tag{36}$$

where $E$ is Young's modulus, and $\mu$ is Poisson's ratio. It follows from (29) and the effective force state concept, the fluid pressure state in (8) can be written as

$$\underline{\mathcal{T}}_w = \omega \mathbf{1} p_w \mathcal{K}^{-1} \underline{\xi}, \tag{37}$$

where $\mathbf{1}$ is the second-order identity tensor.

Similarly, given (24) and (25), the stabilized fluid flow states in the bulk and the fractured space can be written as

$$\underline{\mathcal{Q}} = \omega q_w \mathcal{K}^{-1} \underline{\xi} + \frac{\mathcal{GC}_3}{\omega_0} \omega \underline{\mathcal{R}}_w. \tag{38}$$

$$\underline{\mathcal{Q}}_f = \omega q_f \mathcal{K}^{-1} \underline{\xi} + \frac{\mathcal{GC}_4}{\omega_0} \omega \underline{\mathcal{R}}_f, \tag{39}$$

where $C_3$ and $C_4$ are two material parameters (i.e., micro-conductivities [42]), and $\underline{R}_w$ and $\underline{R}_f$ are the non-uniform fluid pressure state in the bulk and the fractured space, respectively. The latter two terms are defined as

$$\underline{\mathcal{R}}_w = \underline{\Phi} - \widetilde{\nabla \underline{\Phi}} \xi, \tag{40}$$

$$\underline{\mathcal{R}}_f = \underline{\Phi}_f - \widetilde{\nabla \underline{\Phi}}_f \xi. \tag{41}$$

For a three-dimensional case, $C_3$ and $C_4$ are written as

$$\mathcal{C}_3 = \frac{6k_w}{\pi\delta^4}, \tag{42}$$

$$\mathcal{C}_4 = \frac{6k_f}{\pi\delta^4}. \tag{43}$$

In summary, with (29), (30), (37), (38), and (39), the governing equations (8), (16), (12), and (13) are complete. These equations provide a comprehensive framework for modeling the behavior of porous materials, considering both mechanical and fluid-related aspects. Boundary conditions, including natural and essential boundary conditions, are imposed using the boundary layer method [49]. It is worth noting that our approach is highly adaptable, and we have the flexibility to incorporate advanced constitutive models for geomaterials (e.g., [33, 56–58]) into the μPPM paradigm. This capability enables us to explore more complex material behaviors and better represent the mechanical and hydraulic properties of porous media. In what follows, we present the energy-based bond breakage criterion for modeling crack formation in unsaturated porous media. This criterion will provide insights into the initiation and propagation of cracks within the material.

*2.3. Energy-based bond breakage criterion*

In μPPM, the crack can form naturally when sufficient poromechanics bonds break at a material point. In this study, we adopt the energy-based bond breakage criterion to model the bond breakage in unsaturated porous media. In this case, the bond-breakage criterion depends on the maximum deformation energy density in a bond [47]. The energy density in the bond $\boldsymbol{\xi}$ is obtained as

$$\mathcal{W} = \int_0^t \left(\underline{\mathcal{T}} - \underline{\mathcal{T}}'\right) \dot{\widetilde{\mathcal{U}}} dt + \int_0^t \left(\underline{\mathcal{M}} - \underline{\mathcal{M}}'\right) \dot{\underline{\Omega}} dt, \tag{44}$$

where $t$ is the loading time. The critical energy density for a bond can be determined from the critical energy release rate as

$$\mathcal{W}_{cr} = \frac{4\mathcal{G}_{cr}}{\pi \delta^4}, \tag{45}$$

where $G_{cr}$ is the critical energy release rate. The bond breakage is tracked by a parameter $\underline{\varrho}$ that is defined as

$$\underline{\varrho} = \begin{cases} 1 & \text{for } \mathcal{W} < \mathcal{W}_{cr}, \\ 0 & \text{for } \mathcal{W} \geq \mathcal{W}_{cr}. \end{cases} \tag{46}$$

In μPPM, the local damage variable $\phi$ is used to model the inception and propagation of cracks. The local damage variable at a material point is defined as

$$\varphi = 1 - \frac{\int_{\mathcal{H}} \underline{\varrho}\underline{\omega} d\mathcal{V}'}{\int_{\mathcal{H}} \underline{\omega} d\mathcal{V}'}. \tag{47}$$

As in [47], the space between material points $\boldsymbol{x}$ and $\boldsymbol{x}$' is assumed fractured when $\mathcal{W} \geq \mathcal{W}_{cr}$ as well as $\phi \geq \phi_{cr}$ and $\phi' \geq \phi_{cr}$ for both material points. Such material points are defined as fractured points. Moreover, it is assumed that the fractured points are endowed with bulk and fracture fluid pressures. The two fluid pressures are utilized to model unsaturated fluid flow from the bulk to the fractured space (e.g., (14)). The composite relative displacement state $\widetilde{\underline{U}}$ can be decomposed into two parts, of which the one perpendicular to $\underline{\xi}$ represents the crack opening and the one parallel to $\underline{\xi}$ is the crack dislocation. The crack aperture $c$, which is related to the crack opening can be obtained as

$$\underline{c} = \|\widetilde{\underline{\mathcal{U}}}\| \cos \psi - \|\boldsymbol{x}\|, \tag{48}$$

where $\psi$ is the angle between $\underline{\xi}$ and $\widetilde{\underline{Y}}$ [47]. The crack width at fracture point $\boldsymbol{x}$ can be approximated by the averaged bond apertures of all broken bonds as

$$a_f = \frac{\int_{\mathcal{H}} \underline{\omega} \underline{c} d\mathcal{V}'}{\int_{\mathcal{H}} \underline{\omega} d\mathcal{V}'}. \tag{49}$$

*2.4. Numerical implementation*

In our implementation of the μPPM paradigm, we have employed a hybrid Lagrangian-Eulerian meshfree scheme in the spatial domain and an explicit-explicit dual-direction fractional-step algorithm in the temporal domain. This approach allows us to effectively simulate the behavior of porous materials. The porous material body is discretized into a finite number of mixed material points. These points represent both the solid phase and the fluid phase within the material. The solid phase is tracked using a Lagrangian coordinate, which represents the material's deformation, while the fluid phase is tracked using an Eulerian coordinate relative to the solid phase. This relative spatial description helps us account for the movement and flow of fluids within the porous medium. The computational $\mu$PPM paradigm considers three fundamental unknowns: displacement, micro-rotation, and fluid pressures. These variables are crucial for modeling the behavior of the porous material, as they capture both the mechanical and fluid-related aspects of the system. In the temporal domain, we have implemented a dual-direction fractional-step algorithm. This algorithm splits the fully coupled problems into two distinct components, i.e., the unsaturated fluid flow solver and the solid deformation/fracturing solver. The former addresses the fluid flow within the porous medium. It computes fluid pressures and flow rates, taking into account factors like fluid saturation and hydraulic conductivity. The latter focuses on the deformation and fracturing of the solid phase of the porous material. It considers mechanical forces, stresses, and the initiation and propagation of cracks. The dual-direction staggered algorithm allows for parallel computation. Depending on the specific application, either the fluid flow solver or the solid deformation solver can be invoked first. This flexibility is particularly useful for scenarios like hydraulic fracturing, where it may be necessary to address fluid flow before considering solid deformation. For a detailed description of the algorithms and numerical implementation, readers can refer to [25] for further insights and technical details. By combining these elements, our implementation of the $\mu$PPM paradigm provides a robust framework for simulating the behavior of porous materials, capturing both their mechanical response and fluid dynamics in a comprehensive manner.

## 3. Numerical examples

In this section, we present four numerical examples to validate and demonstrate the efficacy of the fully coupled $\mu$PPM paradigm to simulate dynamic cracking and branching in variably saturated deforming porous media. Example 1 deals with hydraulic fracturing in saturated porous media and validating results with the analytical solution. Example 2 concerns the crack branching in dry porous media under high loading rates. Example 3 deals with the single crack branching in unsaturated porous media under dynamic loads. Example 4 studies multiple crack branching in unsaturated porous media under high loading rates. For all examples, the boundary conditions (i.e., essential and natural boundary conditions) are prescribed through the boundary layer method. It is assumed in all numerical examples that the material points on the boundary layers are free to have micro-rotations. The horizon $\delta$ is considered to be equal to the Cosserat length scale $l$.

### 3.1. Example 1: The KGD problem

This study involves the numerical simulation of fluid pressure-induced cracking within a saturated elastic porous medium. The computational results obtained herein are systematically compared with the analytical solutions for the KGD problem presented in [59, 60]. In the context of the KDG problem, we investigate the propagation of a rectilinear crack originating from a line source situated within a homogeneous and isotropic elastic material. Our modeling approach assumes that the fracturing pressure follows the characteristics of a purely viscous fluid, and it considers laminar fluid flow throughout the porous medium. It is important to note that, in the analytical solutions referenced, the medium is idealized as linearly elastic, with no leakage occurring at the fracture surface. Additionally, this study explores the influence of flux rate variations on the phenomenon of crack branching.

Figure 2 illustrates the model setup used in this numerical simulation example. The initial crack has a length of 0.1 m, as depicted in Figure 2. A constant fluid flow rate of $q = 1 \times 10^{-4}$ m²/min is applied at the left end of the crack. Boundary conditions for water pressure are set to zero on the top, bottom, and right boundaries of the model. The right and left boundaries are constrained in horizontal displacement, while the right, top, and bottom boundaries are constrained in vertical displacement. The numerical discretization divides the specimen into a grid of $100 \times 200$ uniform material points, with a grid spacing of $\Delta x = 0.05$ m. A horizon size of $\delta = 0.2$ m is employed, and the time increment is set to $\Delta t = 1.5 \times 10^{-4}$ min. Consistent with [60], the material parameters utilized in the simulation include a bulk modulus of $K = 14.2$ GPa, a shear modulus of $\mu = 11.3$ GPa, a solid density of $\rho_s = 1800$ kg/m³, an initial porosity of $\varphi_0 = 0.19$, a water density of $\rho_w = 1000$ kg/m³, and a hydraulic conductivity of $k_w = 8 \times 10^{-9}$ m/s. For the micropolar material model, the shear modulus is taken as $\mu_c = 5$ GPa, and the micropolar length scale is set at $l = 0.2$ m. We use a stabilization parameter of $G = 0.5$. In this example, the energy-based bond breakage criterion [60] is applied with $\mathbf{G}_{cr} = 100$ N/m.

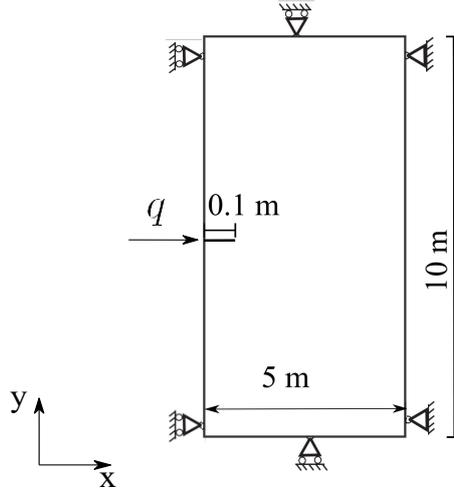

Figure 2: Model setup for example 1.

First, we present the results related to fracturing, fluid flow, and micro-rotation of material points under low flux rate conditions. Figure 3 provides snapshots of the damage variable contour in the deformed configuration at three distinct loading times. A magnification factor of 300 has been applied for clarity. As depicted in Figure 3, the crack exhibits horizontal straight-line growth. Figure 4 displays snapshots of the water pressure contour in a deformed configuration at three loading times. These contours clearly illustrate an increase in water pressure surrounding the growing crack. This phenomenon is attributed to the flow of fluid from the bulk to the fracture space, resulting in a decrease in fracture pressure as the crack lengthens [61]. In Figure 5, snapshots of the micro-rotation of material points at three loading times are plotted. The data presented in Figure 5 reveal that micro-rotation of material points is primarily concentrated at the crack tip, and its magnitude increases with the growth of the crack.

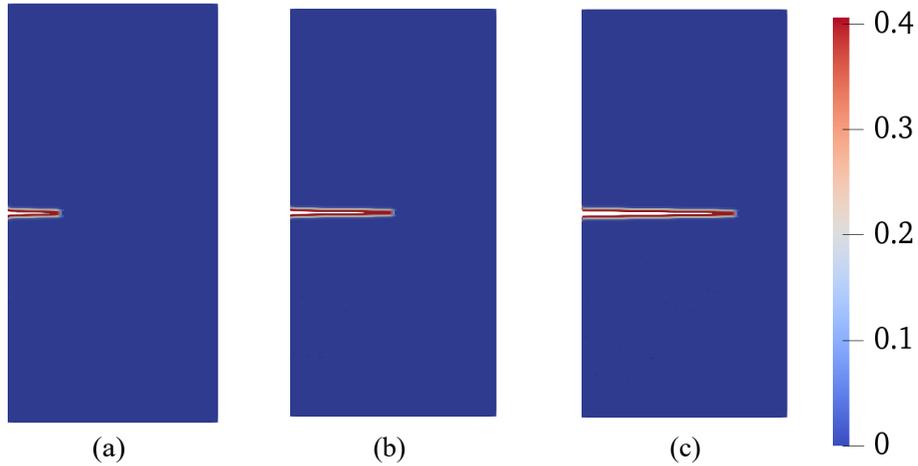

Figure 3: Contours of damage variable in deformed configurations at (a) $t_1 = 4$ min, (b) $t_2 = 7$ min, and (c) $t_3 = 10$ min.

Second, we present a comparison between numerical and analytical results for crack length, crack width, and fracture pressure over the injection period, as detailed in Figures 6, 8, and 7. To facilitate this comparison, we explore two spatial discretization schemes: one with 75×150 points and $\Delta x = 0.067$ m and the other with 100×200 points and $\Delta x = 0.05$ m. Figure 6 showcases the comparison between fracture pressure-time curves obtained through analytical solutions and two discretizations of the $\mu$PPM model. Similarly, Figure 7 illustrates the comparison between fracture width-time curves from analytical solutions and the two discretizations of the $\mu$PPM model. As evident from these results, our numerical simulations align well with the analytical solution. However, it is essential to acknowledge two potential factors contributing to variations in the fracture pressure-time curve. First, the analytical solution assumes a linearly elastic material with no leakage at the fracture interface, whereas the $\mu$PPM model characterizes the material as poroelastic, and the fracture surface allows permeability due to the inclusion of leakage capture within the model. Second, the analytical solution relies on a local theory, while our numerical solution adopts a non-local theory. In the subsequent section, we shift our focus to investigating the impact of a high flux rate on crack branching.

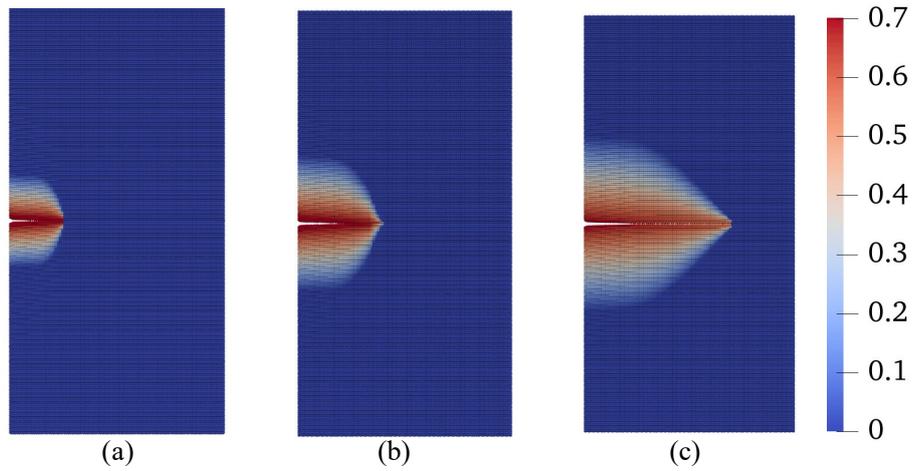

Figure 4: Contours of water pressure (MPa) in deformed configurations at (a) $t_1$ = 4 min, (b) $t_2$ = 7 min, and (c) $t_3$ = 10 min.

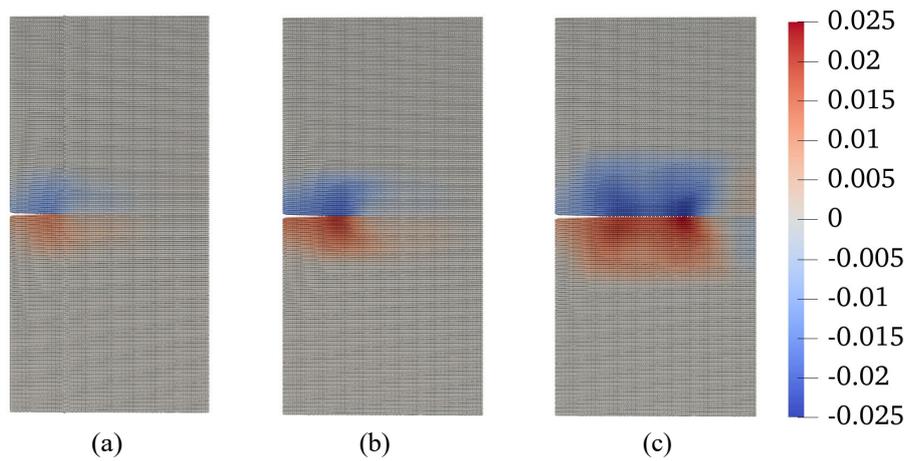

Figure 5: Contours of micro rotation (degree) in deformed configurations at (a) $t_1$ = 4 min, (b) $t_2$ = 7 min, and (c) $t_3$ = 10 min.

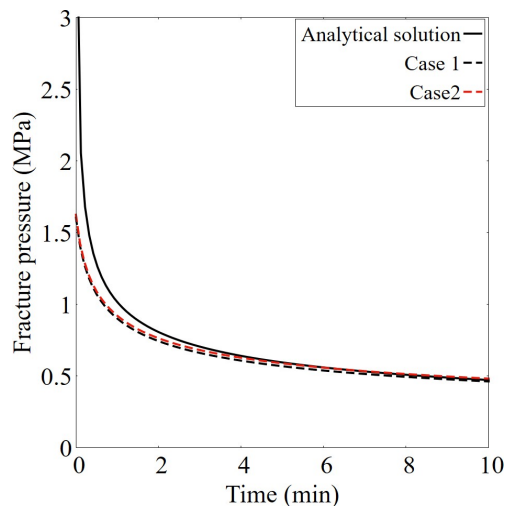

Figure 6: Curves of fracture pressure versus time from the analytical solution and the simulations with two grids.

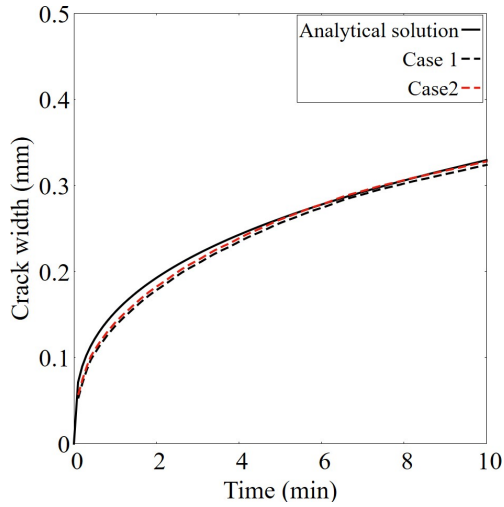

Figure 7: Curves of fracture width versus time from the analytical solution and the simulations with grids.

*3.1.1. Impact of flux rates on the crack branching*

In this part, specifically, we delve into the influence of flux rate on crack branching phenomena in the context of the KGD problem.

First, we investigate the influence of the flux rate on crack branching within this example. The geometry, material parameters, and boundary conditions remain consistent with the KGD problem. Specifically, a fluid flow rate of $q = 2 \times 10^{-3}$ m²/min is applied at the left end of the crack. For this analysis, we employ a spatial discretization scheme with a grid of $75 \times 150$ uniform material points and a grid spacing of $\Delta x = 0.067$ m. Figure 12 provides snapshots of the damage variable contour in the deformed configuration at three distinct loading times. As the flux rate increases, noticeable crack branching phenomena emerge. Figure 13 presents snapshots of water pressure in the deformed configuration at these same three loading times. With an increase in the flux rate, the fluid velocity surges, along with the energy driving the fracture process. When

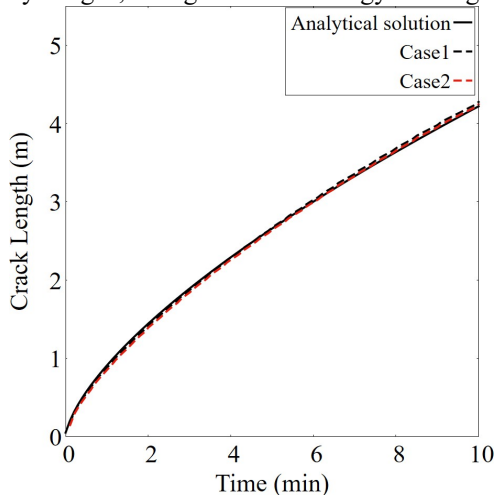

Figure 8: ngth ve        alytical        ations with grids.

the energy       antly h        ipation        rous media, crack
branching        Micr            ial poi        Figure 14, again at
three load       rved in         rotatio        s is predominantly
concentra        k tip a         rack b

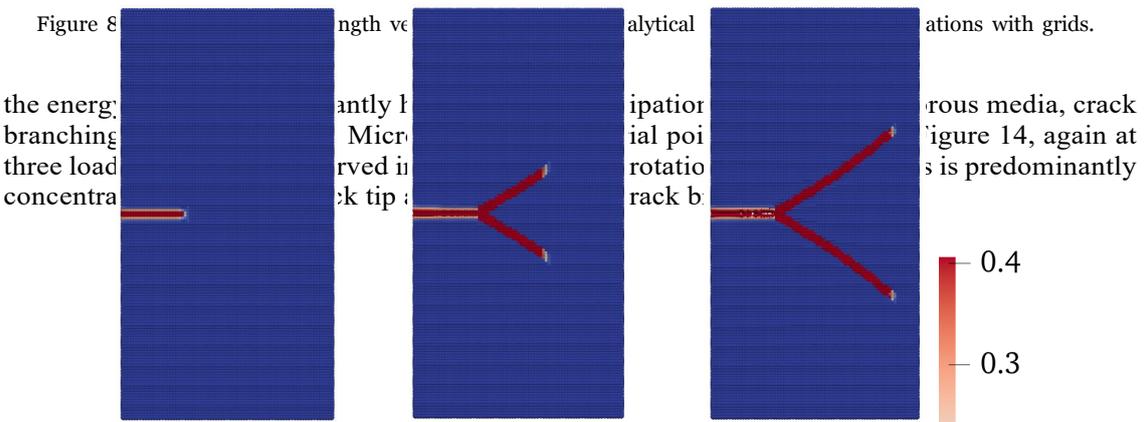

(a) (b) (c)

Figure 9: Contours of crack path in deformed configurations for the case with fluid flow $q = 2 \times 10^{-3}$ m²/min at (a) $t_1 = 3$ min, (b) $t_2 = 4$ min, and (c) $t_3 = 5$ min.

Second, we investigate the impact of spatial discretization on the results while maintaining consistent conditions. To achieve this, we examine two different spatial discretization schemes: one employing 75×150 points with $\Delta x = 0.067$ m and the other utilizing 100 × 200 points with $\Delta x = 0.05$ m. In both cases, the horizon size is set at $\delta = 0.2$ m. Figure 12 presents a comparison of the damage variable contour in the deformed configuration at t=5 min for these two simulation setups. Figure 13 offers a comparison of the water pressure contour in the deformed configuration at t=5 min for both simulations, while Figure 14 compares the micro-rotation of material points in the deformed configuration at the same time instant. Our observations suggest that, under the same horizon size, the choice of spatial discretization scheme has a relatively minor influence on crack branching behavior.

*3.2. Example 2: Crack branching in dry porous media*

This example delves into the phenomenon of crack branching within a dry porous material subjected to high loading rates. Our focus in this study is to investigate how loading rates and spatial discretization impact the occurrence of crack branching. In porous materials, applying high stress within a fracture zone can result in the material's inability to efficiently disperse the energy, leading to the initiation of crack branching. This branching phenomenon typically arises when the crack attains a critical speed of propagation. Furthermore, increasing the loading rate can cause a shift in the failure mode from mode I to a combination of modes. Rapid crack propagation introduces inertia forces at the crack tip, impeding its progression and giving rise to branching. Empirical investigations have demonstrated that crack growth velocity in porous materials tends to be relatively slower. In light of these observations, our analysis centers on understanding the influence of loading rates and spatial discretization on the occurrence and characteristics of crack branching. Additionally, we explore the micro-rotation of material points along the crack path during the branching process. Note that, in this and the subsequent examples (2, 3, and 4), contour variables are presented in deformed configurations, with a magnification factor of 25.

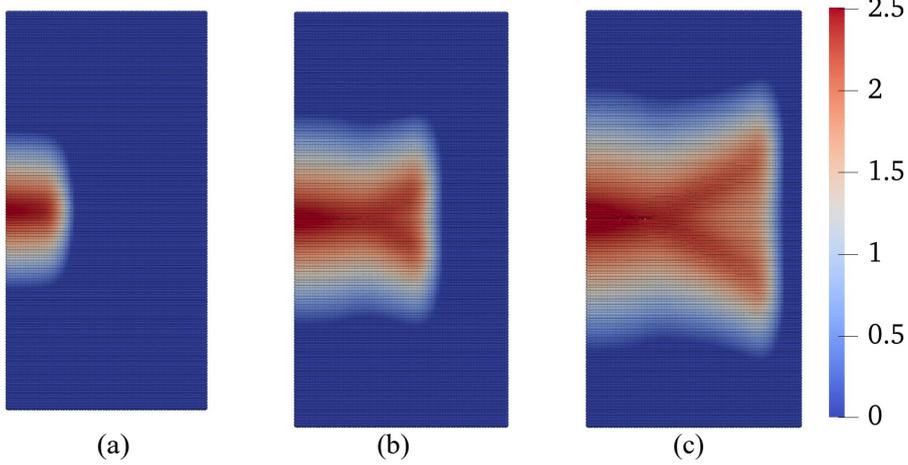

(a) (b) (c)

Figure 10: Contours of water pressure (MPa) for the case with fluid flow $q = 2 \times 10^{-3}$ m²/min at (a) $t_1 = 3$ min, (b) $t_2 = 4$ min, and (c) $t_3 = 5$ min.

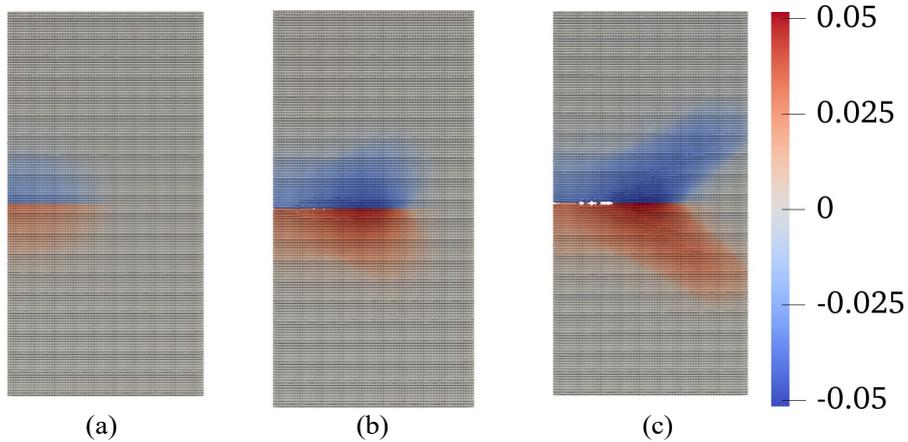

(a) (b) (c)

Figure 11: Contours of micro-rotation (degree) for the case with fluid flow $q = 2 \times 10^{-3}$ m²/min at (a) $t_1 = 3$ min, (b) $t_2 = 4$ min, and (c) $t_3 = 5$ min.

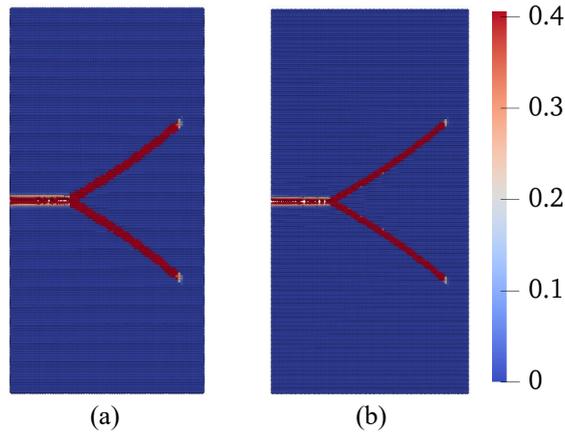

Figure 12: Contours of crack path in deformed configurations at *t* = 5 min from the simulations with (a) grid 1 and (b) grid 2.

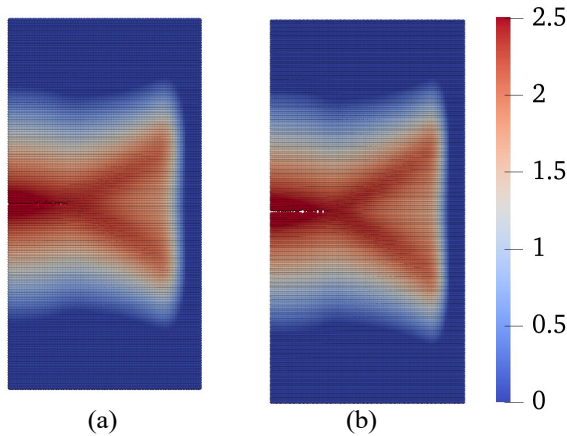

Figure 13: Contours of water pressure (MPa) at *t* = 5 min from the simulations with (a) grid 1 and (b) grid 2

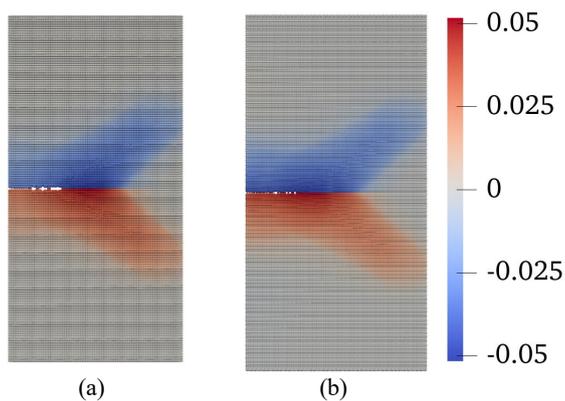

Figure 14: Contours of micro-rotation (degree) at *t* = 5 min from the simulations with (a) grid 1 and (b) grid 2.

Figure 15 illustrates the model setup for this example. The initial crack length is set at 0.5 m, as depicted in Figure 15. At the crack surface, free boundary conditions are applied. Both the left and right boundaries are allowed to deform. Tensile stress in the vertical direction is exerted on the top and bottom boundaries, defined as

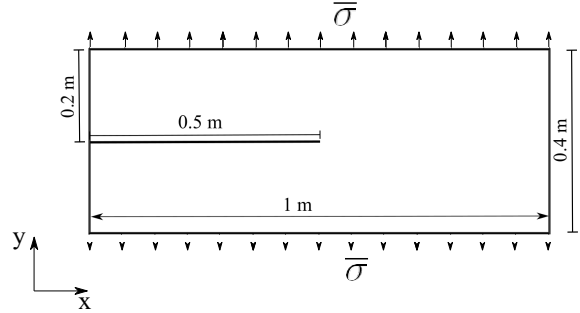

Figure 15: Model setup for example 2.

First, we present the results obtained under a loading rate of $2 \times 10^4$ MPa/s, focusing on the phenomenon of single crack branching. Figure 16 provides snapshots of the damage variable contour in the deformed configuration at three distinct loading stages. As depicted in Figure 16, the loading rate of $2 \times 10^4$ MPa/s induces a single crack branching within the dry porous media. Figure 17 displays snapshots of the micro-rotation of material points in the deformed configuration at the same three loading stages. As illustrated in Figure 17, micro-rotation of material points is primarily concentrated around the crack tip, as well as along the paths of crack propagation and branching. The magnitude of micro-rotation increases with the growth of the crack.

Second, we investigate the impact of a higher loading rate, specifically $4 \times 10^4$ MPa/s, on crack branching. The maximum loading is maintained at $\bar{\sigma}1 = 10$ MPa with $t0 = 2.5 \times 10^{-4}$ s, the same as in the previous case within this example. Figure 18 presents snapshots of the damage variable contour in the deformed configuration at three distinct loading stages. As depicted in Figure 18, the increased loading rate of $4 \times 10^4$ MPa/s results in the occurrence of multiple crack branching events. A notable observation is that increasing the loading rate leads to the initiation of multiple crack branches, as compared to the single crack branching observed in the previous scenario. Figure 19 displays snapshots of the micro-rotation of material points in the deformed configuration at the same three loading stages. Similar to the previous case, micro-rotation of material points is primarily concentrated around the crack tip and along the paths of crack propagation and branching. Furthermore, the magnitude of micro-rotations increases with the growth of the crack, consistent with our earlier observations.

Third, we examine the sensitivity of the numerical results to different spatial discretization schemes within this example. Specifically, we present the outcomes of simulations conducted using two distinct uniform grids. These schemes entail 150×60 points with a grid spacing of $\Delta x = 6.67 \times 10^{-3}$ m (grid 1) and 200×80 points with $\Delta x = 5 \times 10^{-3}$ m (grid 2). In both cases, we maintain the same micro-polar length scale and horizon values, set at 0.02 m. Figure 20 offers a comparison of the damage variable contour in the deformed configuration at $t = 3.5 \times 10^{-3}$ s for the two simulation setups. Likewise, Figure 21 provides a comparison of the micro-rotation of material points in the deformed configuration at the same time instant for both simulations. Notably, our

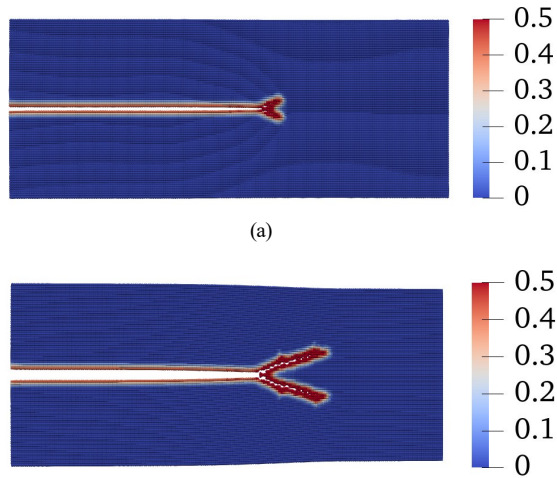

(a)

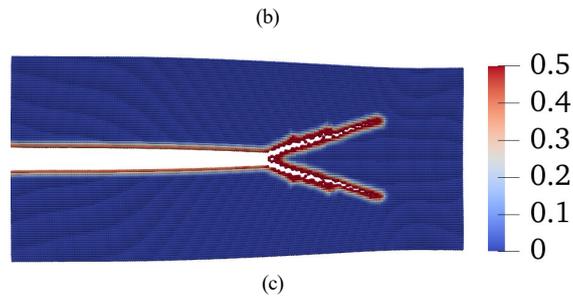
(c)

Figure 16: Contours of damage variable from the simulation with the loading rate $2\times10^4$ MPa/s at (a) $t_1 = 3\times10^{-3}$ s, (b) $t_2 = 3.5 \times 10^{-3}$ s, and (c) $t_3 = 4 \times 10^{-3}$ s.

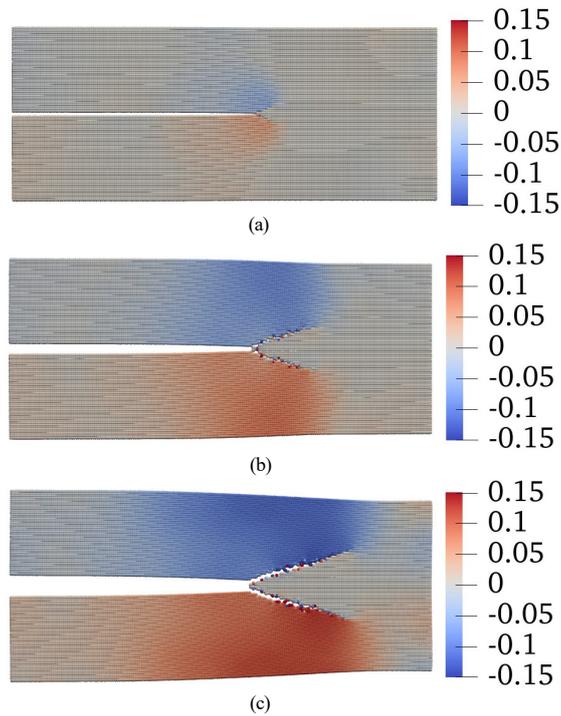

Figure 17: Contours of micro rotation (degree) from the simulation with the loading rate $2 \times 10^4$ MPa/s at (a) $t_1 = 3 \times 10^{-3}$ s, (b) $t_2 = 3.5 \times 10^{-3}$ s, and (c) $t_3 = 4 \times 10^{-3}$ s.

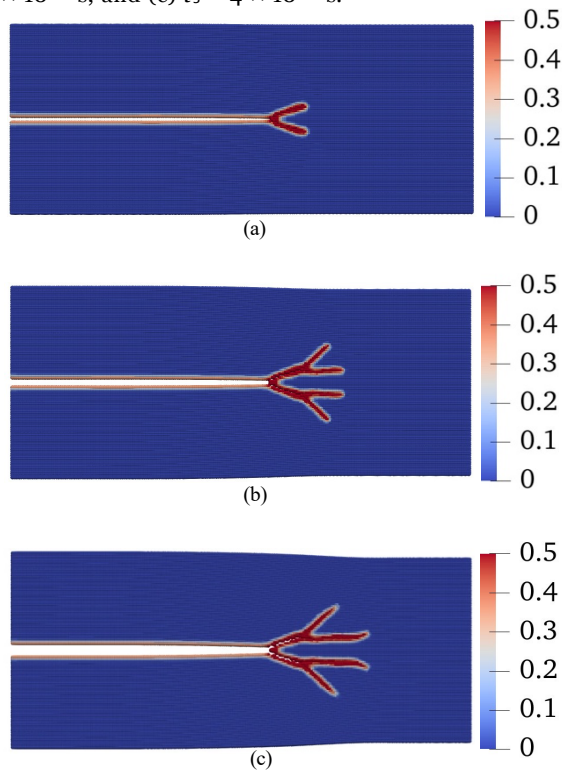

Figure 18: Contours of damage variable from the simulation with the loading rate $4 \times 10^4$ MPa/s at (a) $t_1 = 2.5 \times 10^{-3}$ s, (b) $t_2 = 3.0 \times 10^{-3}$ s, and (c) $t_3 = 3.5 \times 10^{-3}$ s.

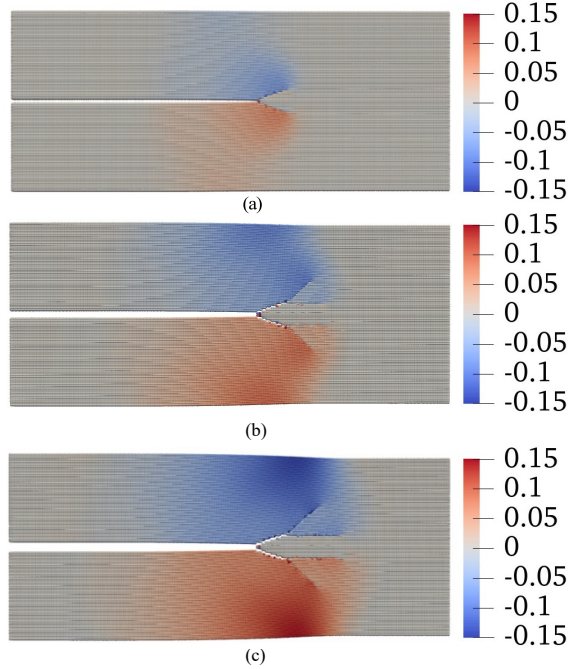

Figure 19: Contours of micro rotation (degree) from the simulation with the loading rate $4 \times 10^4$ MPa/s at (a) $t_1 = 2.5 \times 10^{-3}$ s, (b) $t_2 = 3.0 \times 10^{-3}$ s, (c) $t_3 = 3.5 \times 10^{-3}$ s.

observations from the results suggest that the choice of spatial discretization scheme, under the same horizon size, exerts a relatively minor influence on the occurrence and characteristics of crack branching.

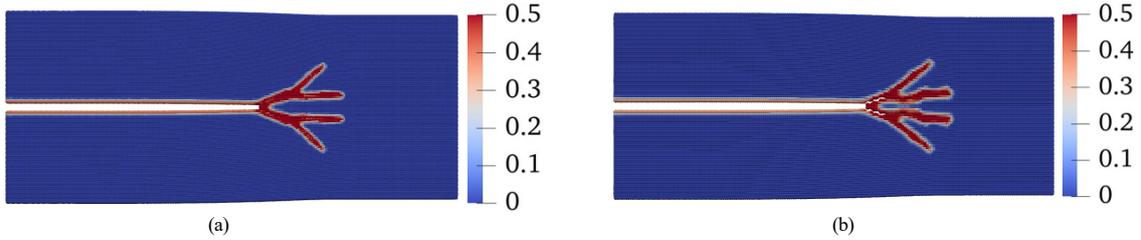

Figure 20: Contours of damage variable from the simulations with the loading rate $2 \times 10^4$ MPa/s at $t = 3.5 \times 10^{-3}$ s: (a) grid 1 and (b) grid 2.

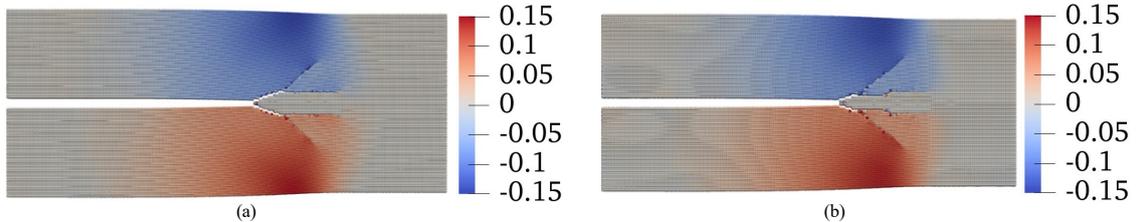

Figure 21: Contours of micro rotation (degree) from the simulations with the loading rate $2 \times 10^4$ MPa/s at $t = 3.5 \times 10^{-3}$ s: (a) grid 1 and (b) grid 2.

### 3.3. Exampl3: Single crack branching in unsaturated porous media

This example focuses on the phenomenon of crack branching within an unsaturated elastic porous material subjected to high loading rates, utilizing the proposed energy-based cracking criterion. Our investigation in this example encompasses the study of fluid flow within unsaturated porous media near the crack tip under conditions of high loading rates. In addition, we explore the influence of spatial discretization on crack branching and fluid flow behavior. The geometry, initial crack, boundary conditions, and material parameters remain consistent with the base simulation in Example 2. Specifically, the maximum loading is set at $\bar{\sigma}_1 = 8$ MPa with $t_0 = 4.0 \times 10^{-4}$ s. For the base simulation, the domain is discretized into $150 \times 75$ material points, employing a grid spacing of $\Delta x = 6.67 \times 10^{-3}$ m. The horizon size is $\delta = 0.02$ m, and a stable time step of $\Delta t = 2.5 \times 10^{-6}$ s is chosen. In the initial state, the water pressure is assumed to be $p_0 = -50$ kPa. To ensure stability in the initial state, the fracture pressure at the initial crack surface is set at $p_{f,0} = -50$ kPa. The water pressure on boundaries is initially -50 kPa, and zero water pressure is gradually imposed on all boundaries. Initially, all boundaries are traction-free, and the total stress is zero, resulting in an effective stress of $\bar{\sigma}_0 = -50$ kPa. The hydraulic parameters for the unsaturated material model include water viscosity $\mu_w = 1 \times 10^{-1}$ Pa·s, water density $\rho_w = 1000$ kg/m³, hydraulic conductivity $k_w = 1 \times 10^{-8}$ m/s, $n = 1.8$, and $s_a = 0.5$ MPa. The stabilization parameter is $G = 0.5$.

First, we present the results obtained from the base simulation, as illustrated in Figures 22, 23, and 24. Figure 22 provides snapshots of the damage variable contour in the deformed configuration at three distinct loading stages. Figure 23 displays snapshots of water pressure in the deformed configuration at the same three loading stages. The contours in Figure 23 reveal that water pressure increases at the branched crack tips while decreasing along the crack path. This behavior is attributed to the flow of water from the bulk space to the fracture space. Initially, water pressure increases and saturates at the crack tip due to local dynamic compression. Subsequently, water pressure flows from the crack tip toward the boundaries, resulting in a decrease in water pressure around the crack tip. Figure 24 illustrates snapshots of micro-rotation of material points at these same three loading stages. As demonstrated in Figure 24, micro-rotation increases as the crack grows.

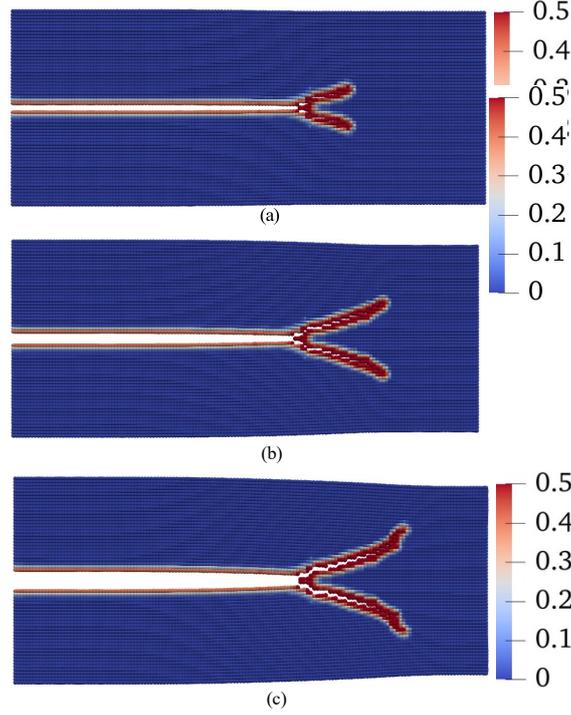

Figure 22: Contours of damage variable from the simulation with the loading rate $2 \times 10^4$ MPa/s at (a) $t_1 = 3.75 \times 10^{-3}$ s, (b) $t_2 = 4 \times 10^{-3}$ s, and (c) $t_3 = 4.25 \times 10^{-3}$ s.

Second, we investigate the influence of spatial discretization on the results while maintaining consistent conditions. To achieve this, we consider two different spatial discretization schemes: one with $150 \times 60$ points and a grid spacing of $\Delta x = 0.067$ m (grid 1), and the other with $200 \times 80$ points and $\Delta x = 0.05$ m (grid 2). In both cases, we assume the same values for the horizon size and micropolar length scale, set at $\delta = 0.2$ m. Figure 25 offers a comparison of the damage variable contour in the deformed configuration at $t = 4 \times 10^{-3}$ s for the two simulation setups. Similarly, Figure 26 provides a comparison of the water pressure contour in the deformed configuration at $t = 4.25 \times 10^{-3}$ s for both simulations, while Figure 27 compares the micro-rotation of material points in the deformed configuration at the same time instant for both simulations. Notably, our observations from the results suggest that, under the same horizon size, the choice of spatial discretization scheme exerts a relatively minor influence on the occurrence and characteristics of crack branching.

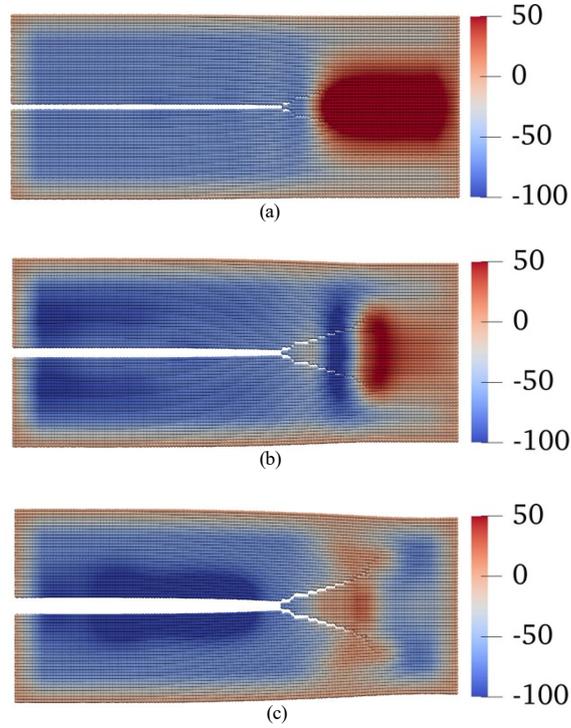

Figure 23: Contours of water pressure (kPa) from the simulation with the loading rate $2 \times 10^4$ MPa/s at (a) $t_1 = 3.75 \times 10^{-3}$ s, (b) $t_2 = 4 \times 10^{-3}$ s, and (c) $t_3 = 4.25 \times 10^{-3}$ s.

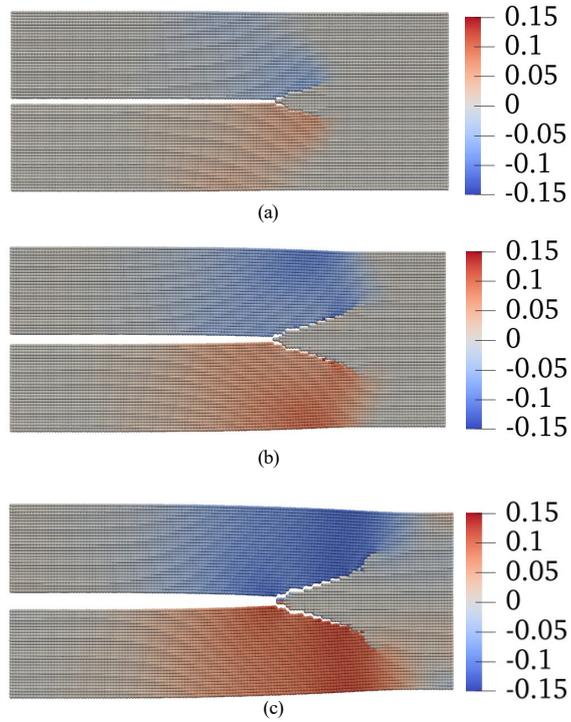

Figure 24: Contours of micro rotation (degree) from the simulation with the loading rate $2 \times 10^4$ MPa/s at (a) $t_1 = 3.75 \times 10^{-3}$ s, (b) $t_2 = 4 \times 10^{-3}$ s, and (c) $t_3 = 4.25 \times 10^{-3}$ s.

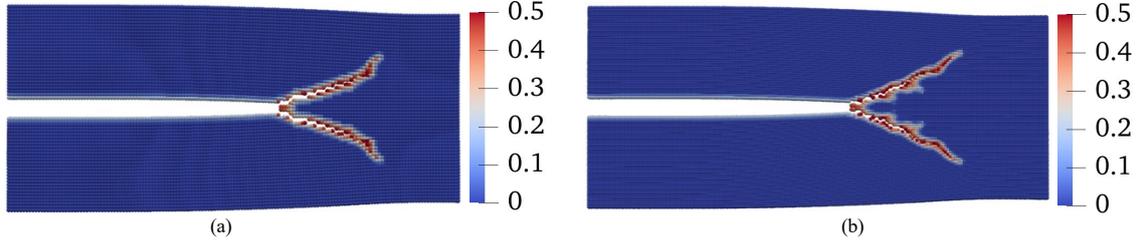

Figure 25: Contours of damage variable from the simulation with the loading rate $2 \times 10^4$ MPa/s at $t = 4.25 \times 10^{-3}$ s: (a) grid 1 and (b) grid 2.

*3.4. Example 4: Multiple crack branching in unsaturated porous media*

This example explores the phenomenon of multiple cracks branching within an unsaturated elastic porous material under high loading rates. The focus of this study includes an examination of fluid flow in unsaturated porous media near the crack tip, with a particular emphasis on understanding the effects of spatial discretization, loading rate, and hydraulic conductivity on multiple crack branching and fluid flow behavior. The geometry, material parameters, and boundary conditions mirror those of the base simulation presented in Example 3. Specifically, the maximum loading is set at $\sigma 1 = 10$ MPa with $t0 = 2 \times 10^{-4}$ s. For the base simulation, the domain is discretized using $200 \times 100$ material points, employing a grid spacing of $\Delta x = 5 \times 10^{-3}$ m. The horizon and micropolar length size remain constant at $\delta = 0.02$ m. A stable time step of $\Delta t = 2.5 \times 10^{-6}$ s is selected.

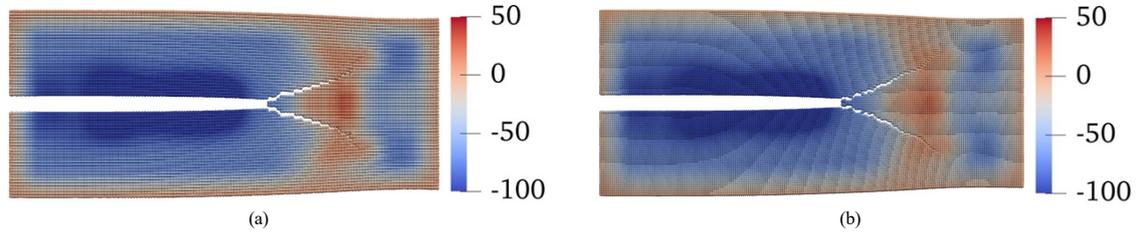

Figure 26: Contours of water pressure (kPa) from the simulation with the loading rate $2 \times 10^4$ MPa/s at $t = 4.25 \times 10^{-3}$ s: (a) grid 1 and (b) grid 2.

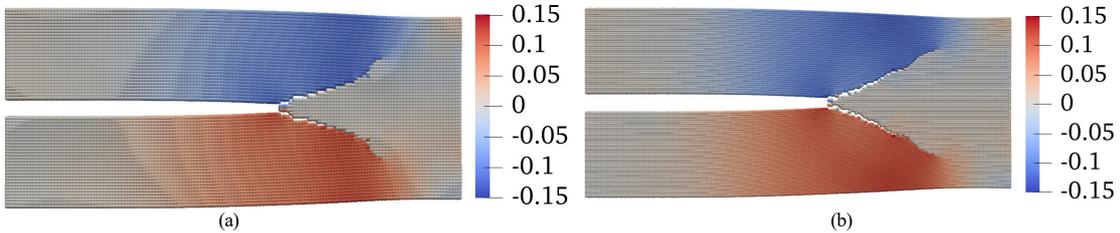

Figure 27: Contours of micro rotation (degree) from the simulation with the loading rate $2 \times 10^4$ MPa/s at $t = 4.25 \times 10^{-3}$ s: (a) grid 1 and (b) grid 2.

First, we present the results obtained from the base simulation, which are illustrated in Figures 28, 29, and 30. Figure 28 provides snapshots of the damage variable contour in the deformed configuration at three distinct loading stages. Figure 29 displays snapshots of water pressure in the deformed configuration at the same three loading stages. The contours in Figure 29 reveal that water pressure increases at the crack tip due to the local volume change rate. Conversely, water pressure decreases along the crack path, attributed to the flow of water from the bulk to the fracture space. Figure 30 illustrates snapshots of micro-rotation of material points at these same three loading stages. As shown in Figure 30, micro-rotation increases as the crack continues to grow.

Next, we explore the impact of spatial discretization on multiple crack branching under consistent conditions. To achieve this, we consider two different spatial discretization schemes: one with 150×60 points and a grid spacing of $\Delta x = 0.067$ m (grid 1), and the other with 200×80 points and $\Delta x = 0.05$ m (grid 2). In both cases, we assume the same value for the horizon, set at $\delta = 0.2$ m. Figure 31 provides a comparison of the damage variable contour in the deformed configuration at $t = 4 \times 10^{-3}$ s for the two simulation setups. As illustrated in Figure 31, the loading rate of $4 \times 10^4$ MPa/s results in the occurrence of multiple crack branching. Figure 32 displays a comparison of the water pressure contour in the deformed configuration at the same time instant for both simulations. Similarly, Figure 33 compares the micro-rotation of material points in the deformed configuration at $t = 4 \times 10^{-3}$ s for both simulations. Notably, our observations from the results suggest that, under the same horizon size, the choice of spatial discretization scheme has a minimal influence on the characteristics of crack branching. In the following sections, we delve into the influence of loading rate and hydraulic conductivity on the phenomenon of multiple crack branching.

*3.4.1. Effect of loading rates on crack branching*

In this section, we investigate the influence of loading rate on the simulation results while maintaining consistent conditions. We consider two loading rates: $5 \times 10^4$ MPa/s and $6 \times 10^4$ MPa/s, with a maximum tensile load set at $\bar{\sigma}_1 = 10$ MPa. All other parameters and boundary conditions remain identical to those of the base simulation in this example. Figure 34 provides a comparison of the damage variable contour in the deformed configuration at time $t = 4 \times 10^{-3}$ s for the two simulations. It is evident from Figure 34 that the crack length is longer when subjected to the higher loading rate. Figure 35 displays a comparison of the water pressure contour in the deformed configuration at the same time instant for both simulations. As demonstrated, in Figure 35, water pressure is concentrated at the crack tip in both simulations. Figure 36 illustrates a comparison of the micro-rotation of material points in the deformed configuration at time $t = 4 \times 10^{-3}$ s for the two simulations. As shown in Figure 36, micro-rotation is concentrated at the crack tip and along the crack paths. Additionally, the magnitude of micro-rotations increases with the growth of the crack. These observations highlight the significant influence of loading rate on crack length and associated fluid flow and micro-rotation behaviors, even under consistent boundary conditions and material parameters.

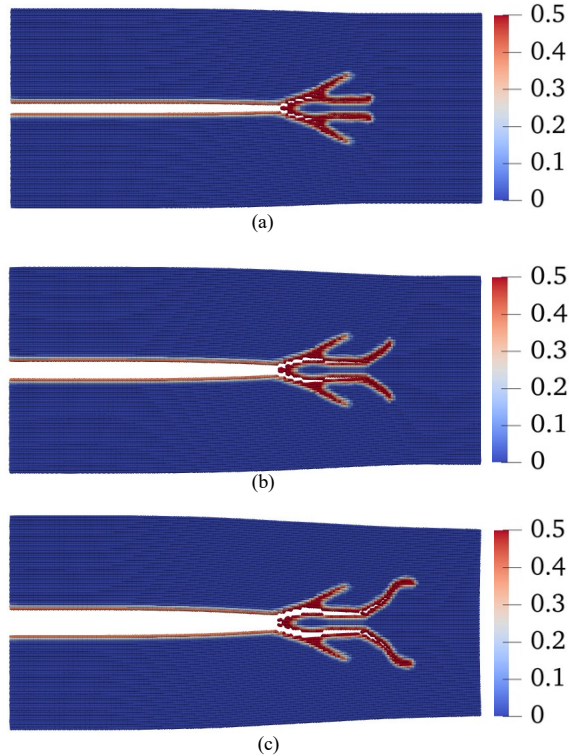

Figure 28: Contours of damage variable from the simulation with the loading rate $5 \times 10^3$ MPa/s at (a) $t_1 = 3.5 \times 10^{-3}$ s, (b) $t_2 = 4 \times 10^{-3}$ s, and (c) $t_3 = 4.5 \times 10^{-3}$ s.

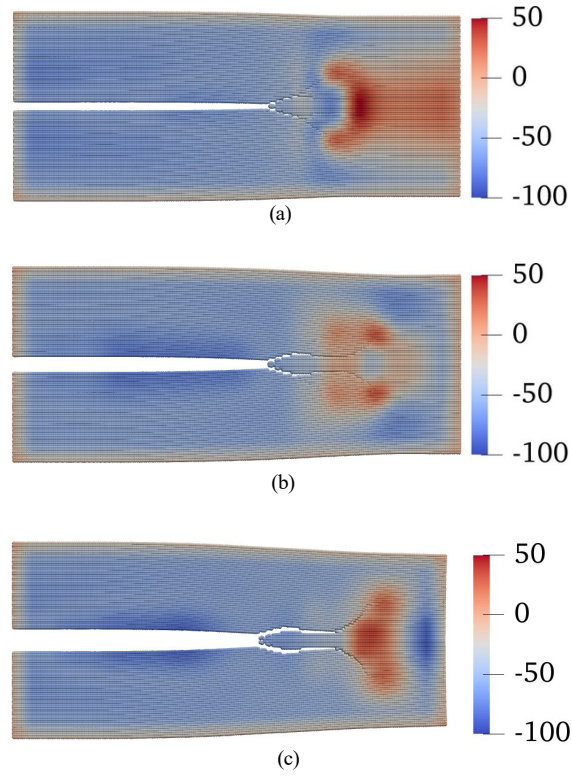

Figure 29: Contours of water pressure (kPa) from the simulation with loading rate $5 \times 10^3$ MPa/s at (a) $t_1 = 3.5 \times 10^{-3}$ s, (b) $t_2 = 4 \times 10^{-3}$ s, and (c) $t_3 = 4.5 \times 10^{-3}$ s.

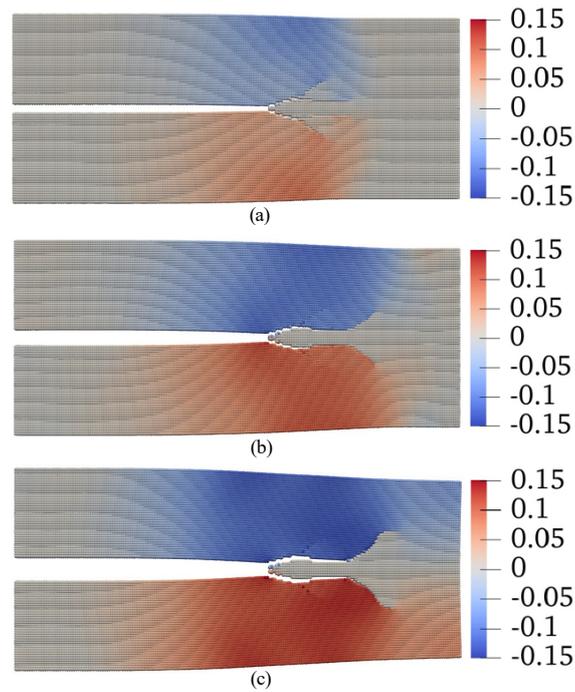

Figure 30: Contours of micro rotation (degree) from the simulation with loading rate $5 \times 10^4$ MPa/s at (a) $t_1 = 3.5 \times 10^{-3}$ s, (b) $t_2 = 4 \times 10^{-3}$ s, and (c) $t_3 = 4.5 \times 10^{-3}$ s.

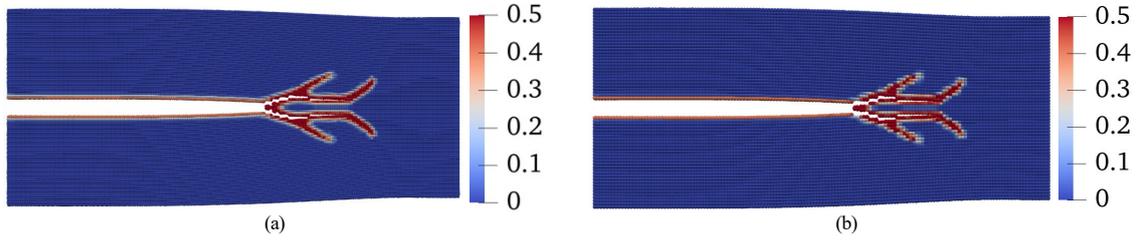

Figure 31: Contours of damage variable from the simulation with loading rate $5 \times 10^4$ MPa/s at $t = 4 \times 10^{-3}$ s: (a) grid 1, and (b) grid 2.

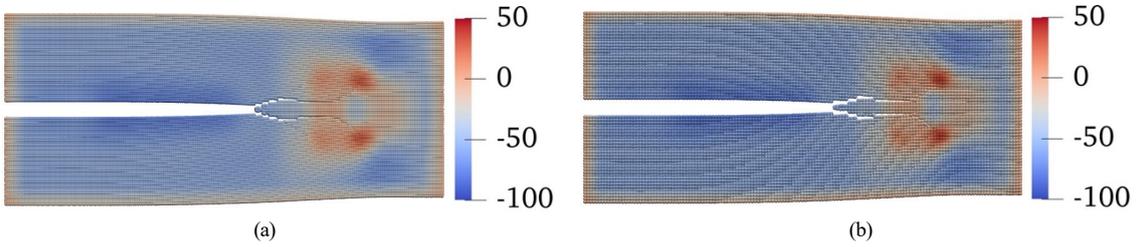

Figure 32: Contours of water pressure (kPa) in from the simulation with loading rate $5 \times 10^4$ MPa/s at $t = 4 \times 10^{-3}$ s: (a) grid 1, and (b) grid 2.

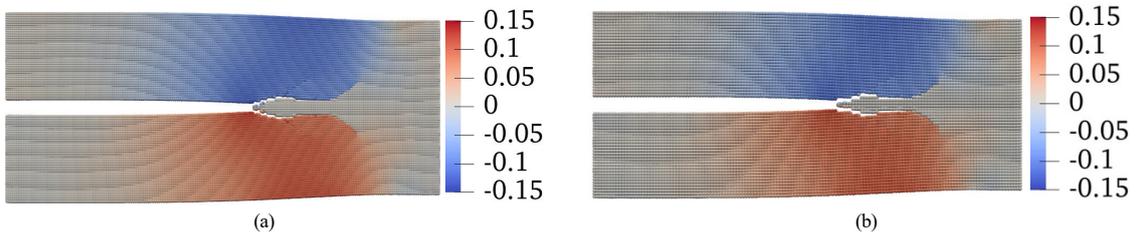

Figure 33: Contours of micro rotation (degree) from the simulation with loading rate $5 \times 10^4$ MPa/s at $t = 4 \times 10^{-3}$ s: (a) grid 1, and (b) grid 2.

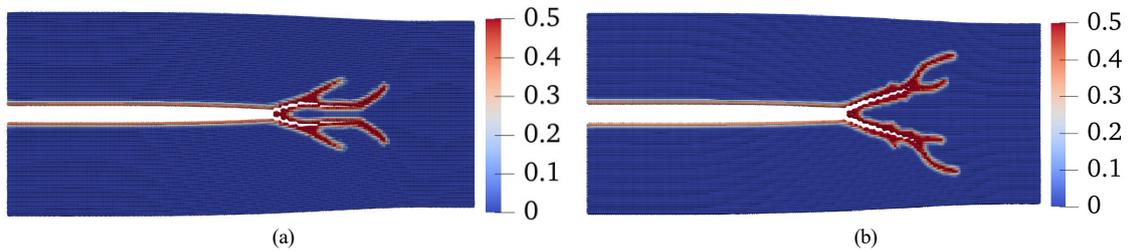

Figure 34: Contours of damage variable at time $t = 4 \times 10^{-3}$ s from the simulations with loading rate: (a) $5 \times 10^4$ MPa/s, and (b) $6 \times 10^4$ MPa/s.

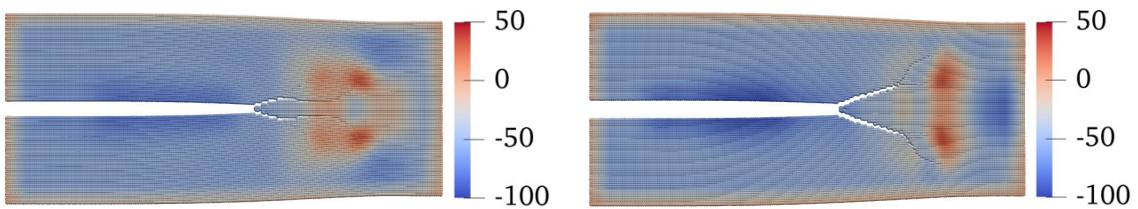

(a) (b)

Figure 35: Contours of water pressure (kPa) at time $t = 4 \times 10^{-3}$ s from the simulations with loading rate: (a) $5 \times 10^4$ MPa/s, and (b) $6 \times 10^4$ MPa/s.

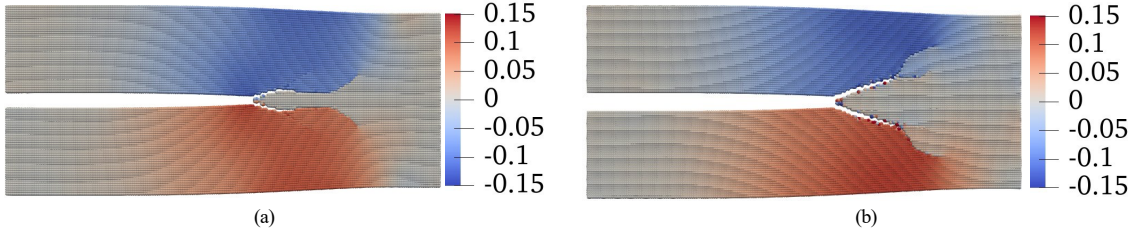

(a) (b)

Figure 36: Contours of micro rotation (degree) at time $t = 4 \times 10^{-3}$ s from the simulations with loading rate: (a) $5 \times 10^4$ MPa/s, and (b) $6 \times 10^4$ MPa/s.

### 3.4.2. Effect of hydraulic conductivity on crack branching

In this section, we investigate the influence of hydraulic conductivity on crack branching, considering three different hydraulic conductivity values: $k_w = 1 \times 10^{-9}$ m/s, $1 \times 10^{-8}$ m/s, and $1 \times 10^{-7}$ m/s. The maximum tensile load is set at $\sigma_1 = 10$ MPa, and the loading rate is maintained at $4 \times 10^4$ MPa/s. The domain is discretized with $200 \times 100$ material points, employing a grid spacing of $\Delta x = 5 \times 10^{-3}$ m. The horizon size remains constant at $\delta = 0.02$ m, and a stable time step of $\Delta t = 2.5 \times 10^{-6}$ s is used. All other parameters and boundary conditions remain consistent with those of the base simulation.

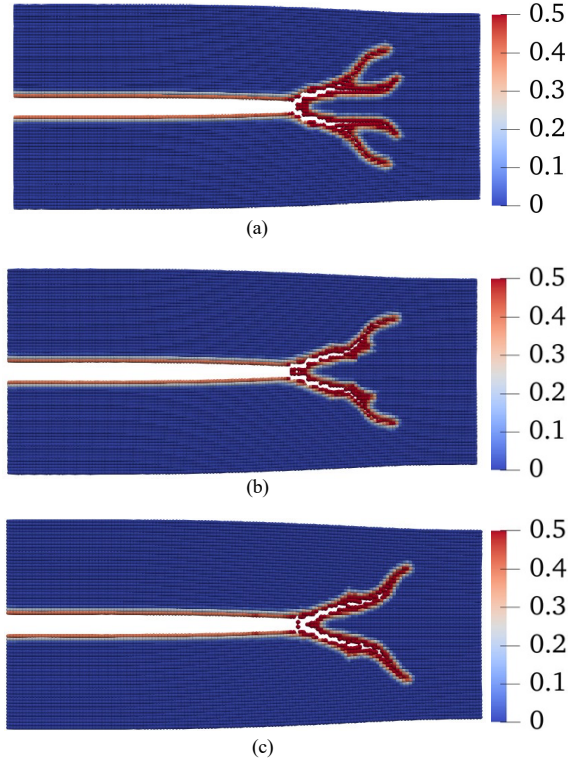

Figure 37: Contours of damage variable in deformed configuration at time $t = 4 \times 10^{-3}$ s from the simulations with (a) $k_w = 1 \times 10^{-9}$ m/s, (b) $k_w = 1 \times 10^{-8}$ m/s, and (c) $k_w = 1 \times 10^{-7}$ m/s.

Figure 37 presents a comparison of the damage variable contour in the deformed configuration at time $t = 4 \times 10^{-3}$ s for the three simulations. As observed in Figure 37, hydraulic conductivity significantly influences the branching behavior. Lower hydraulic conductivity impedes the flow of water pressure from the crack tip to the boundaries, resulting in increased stress concentration at the crack tip and multiple branching events. Figure 38 displays a comparison of the water pressure contour in the deformed configuration at the same time instant for the three simulations. It is evident from Figure 38 that water pressure is concentrated at the crack tip, and increasing hydraulic conductivity leads to higher water pressure values at the crack tip. Figure 39 illustrates a comparison of the micro-rotation of material points in the deformed configuration at time $t = 4 \times 10^{-3}$ s for the three simulations. As shown in Figure 39, micro-rotation is concentrated at the crack tip and along different crack paths for the three simulations. These results emphasize the significant influence of hydraulic conductivity on crack branching behavior, water pressure distribution, and micro-rotation, highlighting the role of fluid flow dynamics in multiple crack branching phenomena.

## 4. Closure

In this paper, we investigate the dynamic branching of cracks in saturated and unsaturated porous media using numerical methods. Our approach builds upon the recently developed coupled micro-periporomechanics ($\mu$PPM) paradigm, an extension of the periporomechanics model. This extended framework accounts for the microrotation of the solid skeleton, enhancing our understanding of the complex interactions within porous media. Within this paradigm, each material point is characterized by three degrees of freedom: displacement, micro-rotation, and fluid pressure. The numerical $\mu$PPM paradigm is mesh-free, harnessing an explicit-explicit split solution algorithm. The coupled $\mu$PPM paradigm incorporates a material length scale based on microstructures, considering micro-rotations of the solid skeleton in alignment with the Cosserat continuum theory for solids. To address the multiphase zero-energy mode instability inherent in the proposed $\mu$PPM, we adopt a stabilized Cosserat $\mu$PPM correspondence principle, which includes unsaturated fluid flow. We introduce an energy-based damage model tailored to the $\mu$PPM paradigm, allowing us to effectively model fracturing in porous media. We present numerical examples to validate and demonstrate the effectiveness of the proposed $\mu$PPM paradigm in modeling fracturing in unsaturated porous media. Our numerical examples explore various scenarios of fluid-driven and deformation-driven crack branching in porous media. To ensure the robustness of our results, we demonstrate the uniqueness of outcomes by employing different discretizations while maintaining the same length scale. Additionally, our numerical examples allow us to analyze the factors influencing crack branching in saturated and unsaturated porous media, including fluid flux rate, loading rates, and hydraulic conductivity. This comprehensive investigation sheds light on the intricate processes of crack branching in porous media and provides valuable insights into the role of various factors in shaping these phenomena.

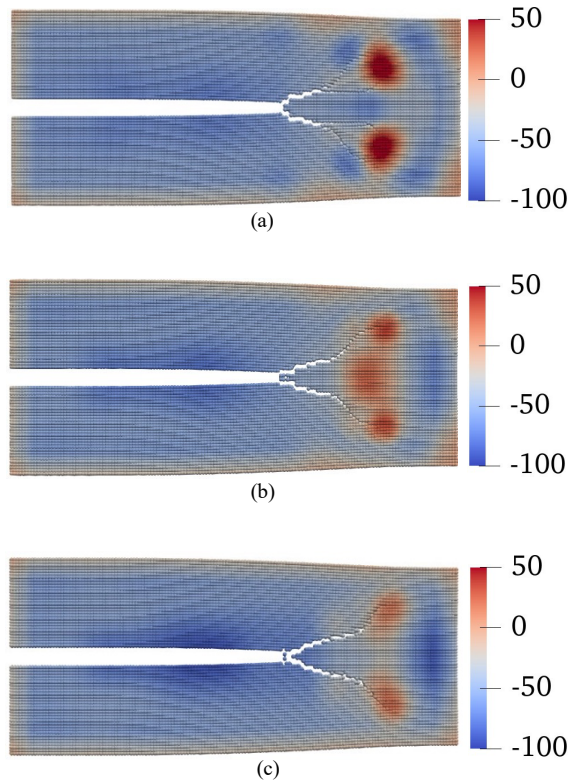

(a)

(b)

(c)

Figure 38: Contours of water pressure (kPa) in deformed configuration at time $t = 4 \times 10^{-3}$ s from the simulations with (a) $k_w = 1 \times 10^{-9}$ m/s, (b) $k_w = 1 \times 10^{-8}$ m/s, and (c) $k_w = 1 \times 10^{-7}$ m/s.

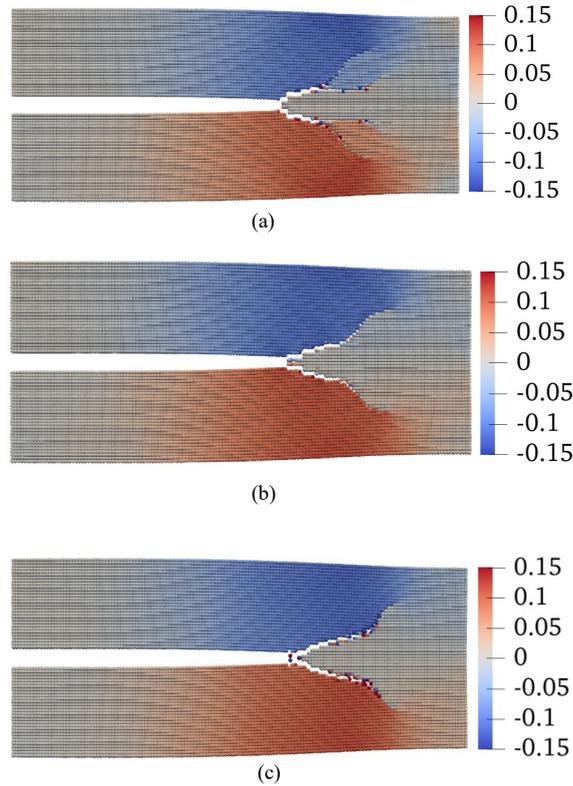

(a)

(b)

(c)

Figure 39: Contours of micro rotation (degree) in deformed configuration at time $t = 4 \times 10^{-3}$ s from the simulations with (a) $k_w = 1 \times 10^{-9}$ m/s, (b) $k_w = 1 \times 10^{-8}$ m/s, and (c) $k_w = 1 \times 10^{-7}$ m/s.

## Acknowledgment

This work has been supported by the US National Science Foundation under contract number 1944009 (CAREER project). The support is gratefully acknowledged. Any opinions or positions expressed in this article are those of the authors only and do not reflect any opinions or positions of the NSF.

## Conflict of Interest Statement

The authors declare no potential conflict of interest.

## Data Availability Statement

The data that support the findings of this study are available from the corresponding author upon reasonable request.